\documentclass[smallcondensed,numbook]{svjour3}

\usepackage[left=3cm, right=3cm, top=1.7cm]{geometry}
\usepackage{amsmath}
\usepackage{amsfonts}
\usepackage{amssymb}
\usepackage[noend]{algpseudocode}
\usepackage{algorithm}
\usepackage{enumitem}
\usepackage{hyperref}

\usepackage{amssymb}
\usepackage{amscd}
\usepackage{epsfig}
\usepackage[table]{xcolor}

\usepackage{graphicx, subfigure}
\usepackage{enumitem}
\usepackage{pgfplots}
\pgfplotsset{width=7cm, height=5cm, compat=1.18}

\usepackage{titlesec}
\titleformat{\chapter}
{\normalfont\LARGE\bfseries}{\thechapter}{1em}{}
\titlespacing*{\chapter}{0pt}{3.5ex plus 1ex minus .2ex}{2.3ex plus .2ex}

\usepackage{multirow}

\newcommand{\textoverline}[1]{$\overline{\mbox{#1}}$}
\newcommand{\ER}[0]{Erd\H{o}s-R\'{e}nyi }

\begin{document}

\title{A Cross Entropy Approach to the Domination Problem and its Variants}
\author{Ryan Burdett, Michael Haythorpe, Alex Newcombe}

\institute{Ryan Burdett
\at Flinders University, 1284 South Road, Tonsley Park, SA, Australia\\
\email{ryan.burdett@flinders.edu.au} \and Michael Haythorpe (Corresponding author)
\at Flinders University, 1284 South Road, Tonsley Park, SA, Australia, Ph: +61 8 8201 2375, Fax: +61 8 8201 2904\\
\email{michael.haythorpe@flinders.edu.au} \and Alex Newcombe
\at Flinders University, 1284 South Road, Tonsley Park, SA, Australia\\
\email{alex.newcombe@flinders.edu.au}} \maketitle

\begin{abstract}The domination problem and several of its variants (total domination, 2-domination and secure domination) are considered. These problems have various real-world applications, but are NP-hard to solve to provable optimality, making fast heuristics for these problems desirable. There is a wealth of highly-developed heuristics and approximation algorithms for the domination problem, however such heuristics are much less common for variants of the domination problem. We redress this by proposing an implementation of the cross entropy method that can be applied to any sensible variant of domination. We present results from experiments which demonstrate that this approach can produce good results in an efficient manner even for larger graphs, and that it works roughly as well for any of the domination variants considered.

\end{abstract}

\section{Introduction}\label{sec-intro}

In this paper we consider simple graphs $G$ which contain a vertex set $V$ and edge set $E$, where $n = |V|$ and $m = |E|$. A subset of the vertices $S \subseteq V$ is a {\em dominating set} for $G$ if every vertex in $V \setminus S$ is adjacent to one or more vertices in $S$. Then, the {\em domination problem} for a graph $G$ is to determine the size of the smallest dominating set in $G$. The latter is denoted as $\gamma(G)$, and is known as the {\em domination number} of $G$.

The domination problem has been studied for the best part of a century, but gained additional interest in the 1970s when the closely related {\em dominating set problem} (for some given constant $K$, determining whether there exists a dominating set with fewer than $K$ vertices) was shown to be NP-complete \cite{garey1979computers}. Hence, the domination problem is NP-hard, and there is no known efficient algorithm which is guaranteed to solve it. Despite this, a number of algorithms have been developed for the domination problem. The best exact algorithms have exponential solving time, such as the $O(1.4969^n)$ algorithm by van Rooij et al. \cite{vanrooij2011exact} which determines not only the domination number, but also the number of minimum dominating sets. The domination problem can also be formulated quite easily as a binary integer programming problem, which software such as CPLEX running on modern desktops are capable of solving even for instances with hundreds of vertices. Beyond exact algorithms, there are several approximation algorithms \cite{mira2022polynomial,parekh1991analysis} and fast heuristics \cite{campan2015fast,casado2023iterated,eubank2004structural}.

In this paper, we add to this literature by adapting the {\em Cross Entropy} (CE) method to the domination problem. CE was first proposed in 1997 by Rubinstein as a means of solving rare event probability estimation, based on the method from \cite{rubinstein1997optimization}. It was subsequently expanded to solve combinatorial optimisation problems \cite{rubinstein1999cross} and has been successfully applied to various NP-hard graph problems in particular, such as traveling salesman problem \cite{deboer2005tutorial} and Hamiltonian cycle problem \cite{eshragh2011hybrid}. However, to the best of the authors' knowledge it has not previously been applied to the domination problem. %TODO might be more but couldn't find very easily

In addition to the (standard) domination problem, we will also show that the general approach we describe in this paper can be adapted to variants of the domination problem. There are many such variants of the domination problem, and we will consider three prominent variants in particular.

The remainder of this paper is laid out as follows. In Section \ref{sec-dom}, we list the variants of the domination problem that we will consider. In Section \ref{sec-ce}, we describe our implementation of the CE method. In Section \ref{sec-results}, we provide experimental results showing the performance of our implementation. Finally, we make some concluding remarks in Section \ref{sec-conclusions}.

\section{Variations of the Domination method}\label{sec-dom}

As described in Section \ref{sec-intro}, a subset $S \subseteq V$ is a dominating set if every vertex in $V \setminus S$ is adjacent to one or more vertices in $S$. Then the domination problem is to find the size of the smallest dominating set in a given graph $G$. We now describe three variants of dominating sets.

Consider a graph $G$ containing vertex set $V$ and edge set $E$. A subset $S \subseteq V$ is a:

\begin{itemize}\item {\bf\em total dominating set} for $G$ if every vertex in $V$ is adjacent to one or more vertices in $S$.
\item {\bf\em 2-dominating set} for $G$ if every vertex in $V \setminus S$ is adjacent to two or more vertices in $S$.
%\item {\bf\em connected dominating set} for $G$ if $S$ is dominating, and the subgraph induced by $S$ is connected.
\item {\bf\em secure dominating set} for $G$ if $S$ is dominating, and for every vertex $v \in V \setminus S$, there exists a vertex $w \in S$ such that $vw \in E$, and $(S \setminus \{w\}) \cup \{v\}$ is a dominating set.
\end{itemize}

For each of the above variants of dominating sets, the corresponding problem is to find the size of the smallest such set. As with the domination problem, some different algorithms have been developed for each of these variants, although the majority of these have been focused on solving particular kinds of graphs. For general graphs, there are few results. Burger et al. \cite{burger2013binary} showed that several variants of the domination problem (including total and secure domination) can be formulated as binary programming problems which can then be solved to optimality using a solver such as CPLEX. It is easy to adapt such formulations to produce a similar binary programming formulation for 2-domination. Foerster \cite{foerster2013approximating} has proposed an approximation algorithm for $k$-domination, while Chleb\'{\i}k and Chleb\'{\i}kov\'{a} \cite{chlebik2008approximation} showed that an approximation algorithm for minimum set cover can be adapted to provide an approximation algorithm for total domination. Unfortunately, the literature on variants of the domination problem is sprawling and often esoteric, typically combining several variants at once to produce results applicable to only very specific situations. This makes it challenging to identify existing algorithms for a desired specific variant. This, in large part, is the motivation for the present work; we seek to propose a framework that can be applied to essentially any variant of domination without requiring that variant to be individually analysed.

For more information on the three variants of domination considered in this paper, we refer the interested reader to Henning and Yeo's excellent 2013 book on total domination \cite{henning2013total} and the 2012 survey on k-domination (and k-independence) by Chellali et al. \cite{chellali2012k}. As of yet there are no books or surveys detailing the state of the literature for secure domination, but a good summary of results is provided in the recent paper by Wang et al. \cite{wang2018complexity}.

\section{Cross Entropy Method}\label{sec-ce}

%ce alg
The premise behind the cross entropy method for discrete optimisation is fairly straightforward. We begin by setting values for some parameters $N, M, \rho, \alpha, r$, the meanings of which will become clear in the following explanation. Then, at iteration $t \geq 0$ a fast heuristic is used to generate $N$ valid solutions according to a probability vector $P^t$, where we define $P^0$ to be a uniform probability vector. The solutions are then ranked by some appropriate scoring function, and we discard all but the best $M$ solutions, which we call the {\em elite set}. A new probability vector $P^*$ corresponding to the elite set is computed; there are various ways this can be done, which we will discuss shortly, and the parameter $\rho$ is used in this calculation to avoid numerical errors. Finally, we compute $P^{t+1} = \alpha P^* + (1 - \alpha) P^t$, and proceed to the next iteration. Throughout this process we keep track of the best solution found so far, and the algorithm terminates once $r$ consecutive iterations have completed without any improvement in the best solution, returning this best solution as the final output. This approach is summarised in Algorithm \ref{ce_alg}.

%The cross entropy method for combinatorial optimisation has a very basic premise, as demonstrated in Algorithm \ref{ce_alg}: we begin with a uniform probability distribution $P_0$ which is used to randomly generate some number $N$ solutions for the problem. We then rank all of the $N$ generated solutions, and choose the best $M$ of them to inform changes to the probability distribution. (TODO explain P*). Finally, we take a linear combination of the previous probability distribution $P_i$ with this new distribution $P^*_i$, to obtain a probability distribution for the next iteration $P_{i+1}$. This linear combination is based on some given value $\alpha \in [0,1]$, in particular $P_{i+1} = (1-\alpha) P_i + \alpha P^*_i$. Throughout each iteration of this process, we keep track of the solution with the highest calculated score. The stopping point of this algorithm is when a given number of iterations $r$ have been completed without improving the score of the best achieved solution.

\begin{algorithm}
\caption{Cross entropy method}\label{ce_alg}
\begin{algorithmic}
    \State Set initial parameters: $N$, $M$, $\rho$, $\alpha$, r
    \State Set initial uniform probablity vector $P^0$
    \State $best \gets \infty$
    \State $t \gets 0$
    \While {$best = \infty$, or $best$ was updated within the past $r$ iterations}
        \State Generate $N$ solutions using $P^t$
        \State Calculate vector $L$ of scores for each solution
        \State Sort solutions and select best $M$ solutions as $elite$ $set$
        \If {$min(L) < best$}
            \State $best \gets min(L)$
        \EndIf
        \State Calculate $P^*$ using $elite$ $set$ and $\rho$
        \State $P^{t+1} \gets \alpha P^* + (1-\alpha)P^t$
        \State $t \gets t+1$
    \EndWhile

    \State {\bf return} $best$
\end{algorithmic}
\end{algorithm}

There are various ways of using the elite set to compute $P^*$, and we briefly summarise the simplest of these now. Suppose that we are using the cross entropy method for a discrete optimisation problem which involves selecting a subset from some universe set satisfying the conditions of the underlying problem. For example, in the context of the domination problem, we need to select a subset of the vertices from the full vertex set so that the selected subset is dominating. Then each entry in $P^*$ corresponds to an element from the universe set, and the purpose of $P^*$ is to give higher probability to entry $i$ if that entry appeared more often in the elite set. Hence, the simplest way of computing $P^*$ is to define $P^*_i$ to be equal to the proportion of solutions from the elite set which contain entry $i$. Afterwards, we normalise $P^*$ to turn it into a probability vector

Of course, there are more sophisticated ways of computing $P^*$. For instance, rather than simply considering the proportion of solutions from the elite set containing an entry $i$, we can take into account the scores of those solutions, affording a stronger probability to entries which appear in the higher quality solutions. It is common to use inverse exponentials for this purpose, and this can sometimes result in some numerical issues on machines if high-precision numbers are not used. The parameter $\rho$ is hence employed in the calculations to avoid the worst of these numerical issues. This is the approach that we will use, as will be explained in the following subsection.

\subsection{Adapting the cross entropy method to the domination problem}

For the domination problem, the universe set in question is the set of vertices in the graph. We can trivially assign a score to any dominating set by setting it equal to the size of the set, with smaller scores being preferable. We will use the notation $L(S) = |S|$ to denote the score of a dominating set $S$.

Recall that at each iteration we will have generated an {\em elite set} of the $M$ best solutions found, that is, the $M$ (out of $N$) dominating sets of smallest size found. We denote this set by $\mathcal{E}$. Furthermore, we will use the notation $\mathcal{E}^i$ to refer to the subset of $\mathcal{E}$ which contains only those dominating sets which include vertex $i$. Then, we will calculate $P^*$ using the following formula for each entry $P^*_i$, where $\delta := \min\limits_{S \in \mathcal{E}}\frac{L(S)}{\log(\rho)}$.

$$P^*_i=\frac{\sum\limits_{\substack{S\in \mathcal{E}^i}}e^{\frac{-L(S)}{\delta}}}{\sum\limits_{S\in \mathcal{E}}e^{\frac{-L(S)}{\delta}}}$$

%$$P^*_i=\left(\sum\limits_{\substack{x\in E\\ x_i=1}}e^{\frac{-L(x)}{\delta}}\right)\left(\sum\limits_{x\in E}e^{\frac{-L(x)}{\delta}}\right)^{-1}$$

Then, all that remains is to address how dominating sets can be generated from a probability vector $P^t$. The process we advocate is as follows. We begin with a set $S$ which is initially empty, and make a copy of $P^t$ which we denote by $P^\dag$. Then we select a vertex by making an observation of a discrete random variable with probability mass function given by $P^\dag$, and add the selected vertex, say $v$, to $S$. Then, we check if $S$ is dominating. If not, we set $P^\dag_v = 0$, normalise $P^\dag$ again, and repeat the process. It is clear that this process will eventually terminate; in the worst case, every single vertex would eventually be added to $S$ which would certainly be a dominating set.

Note that in the above paragraph, we described a process for generating a dominating set, but it is clear that an analagous process could be used for any desired variant of domination, by simply checking the relevant {\em domination criteria} of that variant after each vertex is added. However, it is worth noting that for some variants of domination, the underlying graph may not contain any such sets. For instance, a graph with isolated vertices never contains a total dominating set (in such a case, it is said that the total domination number is $\infty$). For the variants that we consider in this paper, it is clear that every graph contains a dominating set, a 2-dominating set and a secure dominating set. Hence, it is only when considering total dominating sets that we need to first consider the underlying graph.

Although the above approach will always generate a valid dominating set, there is nothing to ensure that the set so generated is minimal. Since it is always beneficial to have minimal sets if possible, we augment the above approach with a second phase which iteratively considers each vertex $v$ in $S$ to see if $S \setminus \{v\}$ also satisfies the relevant domination criteria. If so, that vertex is removed from $S$, and the process continues until all vertices of $S$ have been considered. It is clear that this process results in a minimal dominating set. The only question is in which order we should consider the vertices for removal; note that seeking to do so optimally is an NP-hard problem\footnote{To see this, note that $S$ could contain every vertex at the conclusion of phase 1, in which case solving the second phase optimally is equivalent to solving the domination problem.}. Instead, we propose the following fast heuristic. We again make a copy of $P^t$ and denote it $P^\dag$. Then we define $\overline{P}$ to be the vector containing the following entries:

$$\overline{P}_i := \left\{\begin{array}{ll}1 - P^\dag_i, & \mbox{\;\;\;\;\;\;\;\; if } i \in S,\\
0 & \mbox{\;\;\;\;\;\;\;\; otherwise.}\end{array}\right.$$

After normalising $\overline{P}$, we then use it in the same manner as we previously used $P^\dag$ to randomly generate a sequence of the vertices in $S$, and this is the order used to consider the vertices for removal. We summarise both phases in Algorithm \ref{dom_alg}, where the {\em rand} function corresponds to an observation of a discrete random variable with probability mass function given by the associated probability vector.

\begin{algorithm}
\caption{Algorithm for generating a minimal dominating set given $P^t$}\label{dom_alg}
\begin{algorithmic}
\State Phase 1:
\State $S \gets \emptyset$
\State $P^\dag \gets P^t$
\While{$S$ does not satisfy the relevant domination criteria}
    \State Normalise $P^\dag$.
    \State $v \gets rand(P^\dag)$
    \State $S \gets S \cup \{v\}$
    \State $P^\dag_v \gets 0$
\EndWhile\\
\State Phase 2:
\State $\overline{P} \gets (1-P^t)$
\For{all vertices $v \not\in S$}
    \State $\overline{P}_v \gets 0$
\EndFor
\While{$\overline{P}$ is not a zero vector}
    \State Normalise $\overline{P}$
    \State $v \gets rand(\overline{P})$
    \If{$S \setminus \{v\}$ satisfies the relevant domination criteria}
        \State $S \gets S \setminus \{v\}$
    \EndIf
    \State $\overline{P}_v \gets 0$
\EndWhile\\
\State Return $S$
\end{algorithmic}
\end{algorithm}

%dom alg
% set x empty vector of length n
% copy p_t into pbar
% while domination constraints not met
%   select an index k using probability dist pbar
%   set x_k = 1
%   set pbar_k = 0
% set pbar = 1-p_t
% for t = 0..n-1
%   select an index k using probability dist pbar
%   set pbar_k = 0
%   if x_k not 0
%       set x_k = 0
%       if domination constraints no longer met
%           set x_k = 1

\subsection{Checking the relevant domination criteria}

Note that for the cross entropy implementation described in the previous subsections, we need to check the relevant domination criteria $O(n)$ times for each of the $N$ dominating sets generated per iteration. The standard method of checking if a set is dominating takes $O(m)$ time, and so using this approach would require us to spend $O(nmN)$ time each iteration generating the dominating sets. This is the most computationally expensive component of the algorithm, and so we take steps to perform these checks more efficiently. Specifically, we utilise efficient updating procedures whenever we add or remove a vertex from $S$ that allow us to track how close we are to meeting (or failing) the relevant domination criteria.

For (standard) domination, this is uncomplicated. For each vertex, we can keep track of how many vertices from its closed neighbourhood are in $S$. Whenever a vertex $v$ is added to $S$ (or removed from $S$), we need to update only those vertices adjacent to $v$. If we denote the degree of $v$ to be $d(v)$, this updating procedure will occur in $O(d(v))$ time. Meanwhile, we can maintain a count of undominated vertices that is decreased whenever a vertex becomes dominated (or vice versa). Checking if $S$ is a dominating set is then as simple as checking if this count is equal to zero. This approach requires us to spend $O(mN)$ time per iteration, which represents an order of magnitude improvement in computation time. For 2-domination and total domination, checks analogous to the above are straightforward and provide similar improvement in computation time.

Of the variants of domination considered in this paper, it is only secure domination that requires more thought. The naive method of checking if $S$ is secure dominating involves looking at each vertex $v$ not in $S$, considering each of its neighbours $w$ which are in $S$, and checking if $(S \cup \{u\}) \setminus \{w\}$ is dominating, which requires $O(m^2)$ time. There are more sophisticated methods of checking secure domination; Burger et al. \cite{burger2013two} give one such method, however they do not indicate the computational complexity other than to state it is more complex than checking standard domination. Regardless, we again seek to improve on this by proposing an updating procedure. Like before, for each vertex, we keep track of how many vertices in its closed neighbourhood are in $S$; for vertex $v$, call this number $domcount(v)$. Then, whenever a vertex $v$ is added or removed from $S$, we only need to update vertices adjacent to $v$. Furthermore, we say that $w$ is {\em capable of defending} $v$ if $w \in S \cap N(v)$ and every neighbour of $w$ is dominated in $(S \cup \{v\}) \setminus \{w\}$, and note that this criteria can be efficiently checked as follows. If there is any neighbour $u$ of $w$ such that $domcount(u) = 1$ and $u \not\in N[v]$ then $w$ is not capable of defending $v$; otherwise, it is capable of defending $v$. Then, for each vertex, we keep track of how many vertices in its open neighbourhood are capable of defending it. Whenever a vertex $v$ is added to $S$ (or removed from $S$), we need to update only those vertices within distance 3 of $v$. Note that we can compute which vertices are within distance 3 of each vertex in advance (i.e. before the first iteration of cross entropy begins) so this need not contribute to the runtime of the updating procedure. Finally, we can maintain a count of vertices which are not in $S$ and have no vertices capable of defending them, and then $S$ is a secure dominating set if and only if this count is equal to zero. For graph families where the diameter grows with $n$, this approach represents an order of magnitude improvement in computation time over the course of generating the secure dominating set.

%There are many approaches to this, for example one could represent domination as an integer program and find solutions to this, akin to what was done in (TODO there was this paper that i found that solved CE by representing the problem as an LP of some kind, but i need to go back and find it). The approach we instead decided to use was a dynamic programming method that keeps track of which vertices were dominated at any given stage of generating a dominating set. To put it into more detail, the dominating set generation goes as follows (see Algorithm \ref{dom_alg}): first we start with empty vectors $x$ and $y$ - both of length $n$, where $y_k=1$ means that vertex $k$ is dominated - and a given probability distribution $P_i$ for the current iteration. We make a copy of $P_i$, $\overline{P}$ to be modified in this algorithm. One by one we select a vertex based on the probabilities in $\overline{P}$. When a vertex $k$ is selected we set $x_k=1$ and $\overline{P}_k=0$. We also increment $y_j$ for all vertices $j \in N[k]$. Vertices are repeated added using the modified $\overline{P}$ until $y_k > 0$ for all vertices $k$ in $G$, in which case we know that a set $S$ obtained from $x$ dominates $G$.

\section{Experimental Results}\label{sec-results}
%redefine specific parameters in a list

The implementation of the cross entropy method for domination described in Section \ref{sec-ce} requires the user to input many parameters, and the overall performance (both in terms of solution quality, and computation time) are impacted by these. In order to determine which combinations of parameters to use for the main experiments, initial experiments were conducted to identify good default values. The results of those experiments indicated that the following values provide a good balance between solution quality and computation time for our implementation: $N = 100, M =10, \alpha = 0.2, \rho = 0.01, r = 20$. In particular, we found that increasing the amount of computation performed (i.e. by increasing the number of dominating sets generated per iteration, or by increasing the number of iterations without improvement before stopping) from these default settings offered, at best, only marginal improvements in solution quality.
%\begin{itemize}
%    \item $N$ - total number of dominating sets to generate, per iteration.
%    \item $M$ - number of elite sets
%    \item $\alpha$ - smoothing parameter for updating $P$
%    \item $\rho$ - rarity parameter, used in calculating $P^*$
%    \item $r$ - number of iterations without improvement before terminating
%\end{itemize}
%
%Our experiments indicated that the following values provide good results for this implementation of cross entropy: $N = 100, M = N/10=10, \alpha = 0.2, \rho = 0.01, r = 20$.

%maybe not?? - describe how we did initial tests to determine which parameters to use for more detailed tests

%describe the more detailed tests and show results

Using the above parameters settings, we now present experimental results for our implementation of the cross entropy method on several kinds of graphs. In particular, we consider grid graphs, flower snarks, unit disk graphs, \ER random graphs, and various graphs from the literature. The former two graph families have been chosen since the domination numbers (and some variants thereof) are known, while the remaining graphs have been chosen as they have previously been used to evaluate graph algorithms. Where possible, we obtain results for all four types of domination considered in this paper (domination, 2-domination, total domination and secure domination.) The only exception is for instances with isolated vertices, in which case we do not obtain results for total domination. We now briefly discuss each type of graph we consider.

Square grid graphs $G(n,n)$ are the Cartesian products of two paths with length $n$. They have been considered extensively in the context of domination, with the domination numbers now known for all cases \cite{gonccalves2011domination}. Fascinatingly, the formula for the domination number of rectangular grid graphs $G(m,n)$ contains 23 special cases before settling into a standard formula for $m,n \geq 16$. Values for the other variants of domination are also known in some cases \cite{burdett2020improved,rao20192,soltankhah2010results}.

Flower snarks $J(k)$ are 3--regular graphs containing $4k$ vertices introduced by Isaacs \cite{isaacs1975infinite}. In a recent paper \cite{burdett2023variants} the domination, 2-domination, total domination and secure domination numbers (among others) were determined.

Unit disk graphs \cite{clark1990unit} are generated in the following way. Given parameters $c,r,m,n$, a set of $c$ points are generated uniformly by random within an $m \times n$ rectangle. Then, a graph is produced where the vertices correspond to the points, and an edge exists between two vertices if the distance between their corresponding points is no more than $2r$. Equivalently, one can draw circles of radius $r$ around each point, and then an edge exists between two vertices if their corresponding circles overlap. Of course, it is not guaranteed that such graphs will be connected; in the upcoming Table \ref{udg_table}, we label the graphs as ``UDG\_c-r-m-n\_s", where s corresponds to a random seed, and we only include graphs in our experiments which are connected.

\ER random graphs are generated in the following way. Given parameters $N,p$, a graph with $N$ vertices is generated in which each edge exists with probability $p$. In the interest of studying relatively sparse graphs, we have experimentally chosen values of $p$ to ensure a low average degree for a given $N$, and then generated many graphs with these parameter settings until a graph with average degree very close to the desired value is obtained. As such, in the upcoming Table \ref{random_table}, we label these graphs ``randomN\_d" where d is the desired degree. Note that since these graphs have low average degree, as $N$ increases it is almost certain that there will be some isolated vertices. As such, we do not consider total domination for these graphs.

Both unit disk graphs and \ER random graphs have been considered in many graph contexts, including as experimental instances for domination algorithms. For example, see \cite{casado2023iterated,chalupa2018order} in which domination algorithms are presented. Those papers also consider a number other graphs from the literature, and so we include many of them in our experiments as well, omitting only those which are too large to be computationally feasible. We provide individual citations for the chosen instances in the upcoming Table \ref{dimacs_table}.

For each instance considered, we run the cross entropy method ten times (with ten different random seeds) and then return the best solution produced. Where possible, we will compare the best solutions produced by the cross entropy method to the optimal solutions for those instances. For instances where the optimal solutions are not known from the literature, we use CPLEX to solve binary programming formulations of the relevant domination variant. In particular, we use the formulations for domination and total domination from \cite{burger2013binary}, and the formulation for secure domination from \cite{burdett2020improved}. We also use a formulation for 2-domination analagous to the formulation for domination in \cite{burger2013binary}. We set a time limit of 10,000 seconds for CPLEX to terminate, and if CPLEX is unable to produce an optimal solution by this time, we take the best solution (upper bound) produced to this stage. The experiments were conducted on an Intel(R) Core(TM) i5-12500 CPU with a 6 core, 3.00GHz processor and 16GB of RAM, running Windows 10 Enterprise version 22H2, using a C++ implementation of the cross entropy algorithm. The solutions and upper bounds for the formulations were obtained on the same PC using CPLEX v22.1.0.

\begin{figure}[h!]
    \centering
    {\bf Flower Snarks}\\
    \begin{tikzpicture}
\begin{axis}[
title={Domination}, xlabel={k}, ylabel={$|S|$}, xmin=0, xmax=100,ymin=0, ymax=150, xtick={0, 25, 50, 75, 100}, ytick={0, 50, 100, 150}, ymode=linear, grid style=dashed,
]

\addplot[dotted]
coordinates {
(5,5)(10,10)(15,16)(20,22)(25,27)(30,33)(35,40)(40,46)(45,50)(50,57)(55,64)(60,71)(65,76)(70,82)(75,88)(80,94)(85,102)(90,108)(95,114)(100,119)
};

\addplot[solid]
coordinates {
(5,5)(10,10)(15,15)(20,20)(25,25)(30,30)(35,35)(40,40)(45,45)(50,50)(55,55)(60,60)(65,65)(70,70)(75,75)(80,80)(85,85)(90,90)(95,95)(100,100)
};

\end{axis}
\end{tikzpicture}
    \hskip 5pt
    \begin{tikzpicture}
\begin{axis}[
title={Total Domination}, xlabel={k}, ylabel={$|S|$}, xmin=0, xmax=100,ymin=0, ymax=175, xtick={0, 25, 50, 75, 100}, ytick={0, 50, 100, 150,175}, ymode=linear, grid style=dashed,
]

\addplot[dotted]
coordinates {
(5,8)(10,16)(15,23)(20,32)(25,41)(30,49)(35,58)(40,66)(45,75)(50,85)(55,92)(60,101)(65,111)(70,120)(75,128)(80,138)(85,145)(90,156)(95,165)(100,173)
};

\addplot[solid]
coordinates {
(5,8)(10,16)(15,23)(20,30)(25,38)(30,46)(35,53)(40,60)(45,68)(50,76)(55,83)(60,90)(65,98)(70,106)(75,113)(80,120)(85,128)(90,136)(95,143)(100,150)
};

\end{axis}
\end{tikzpicture}
    \vspace*{0.5cm}
\begin{tikzpicture}
\begin{axis}[
title={2-Domination}, xlabel={k}, ylabel={$|S|$}, xmin=0, xmax=100,ymin=0, ymax=200, xtick={0, 25, 50, 75, 100}, ytick={0, 50, 100, 150, 200}, ymode=linear, grid style=dashed,
]

\addplot[dotted]
coordinates {
(5,9)(10,18)(15,27)(20,36)(25,46)(30,55)(35,64)(40,75)(45,83)(50,94)(55,104)(60,113)(65,124)(70,133)(75,142)(80,153)(85,164)(90,173)(95,183)(100,194)
};

\addplot[solid]
coordinates {
(5,9)(10,18)(15,25)(20,34)(25,43)(30,50)(35,59)(40,68)(45,75)(50,84)(55,93)(60,100)(65,109)(70,118)(75,125)(80,134)(85,143)(90,150)(95,159)(100,168)
};

\end{axis}
\end{tikzpicture}
    \hskip 5pt
\begin{tikzpicture}
\begin{axis}[
title={Secure Domination}, xlabel={k}, ylabel={$|S|$}, xmin=0, xmax=100,ymin=0, ymax=175, xtick={0, 25, 50, 75, 100}, ytick={0, 50, 100, 150, 175}, ymode=linear, grid style=dashed,
]

\addplot[dotted]
coordinates {
(5,8)(10,16)(15,24)(20,33)(25,41)(30,50)(35,59)(40,67)(45,76)(50,85)(55,94)(60,103)(65,112)(70,121)(75,129)(80,139)(85,148)(90,156)(95,165)(100,175)
};

\addplot[solid]
coordinates {
(5,8.0)(10,16.0)(15,23.0)(20,31.0)(25,38.0)(30,46.0)(35,53.0)(40,61.0)(45,68.0)(50,76.0)(55,83.0)(60,91.0)(65,98.0)(70,106.0)(75,113.0)(80,121.0)(85,128.0)(90,136.0)(95,143.0)(100,151.0)
};

\end{axis}
\end{tikzpicture}
    \caption{Results for Flower Snarks $J(k)$. The size of the best solutions returned by the cross entropy method are displayed using a dotted line, while the known optimal values are displayed using a solid line.}
    \label{flower_plots}
\end{figure}
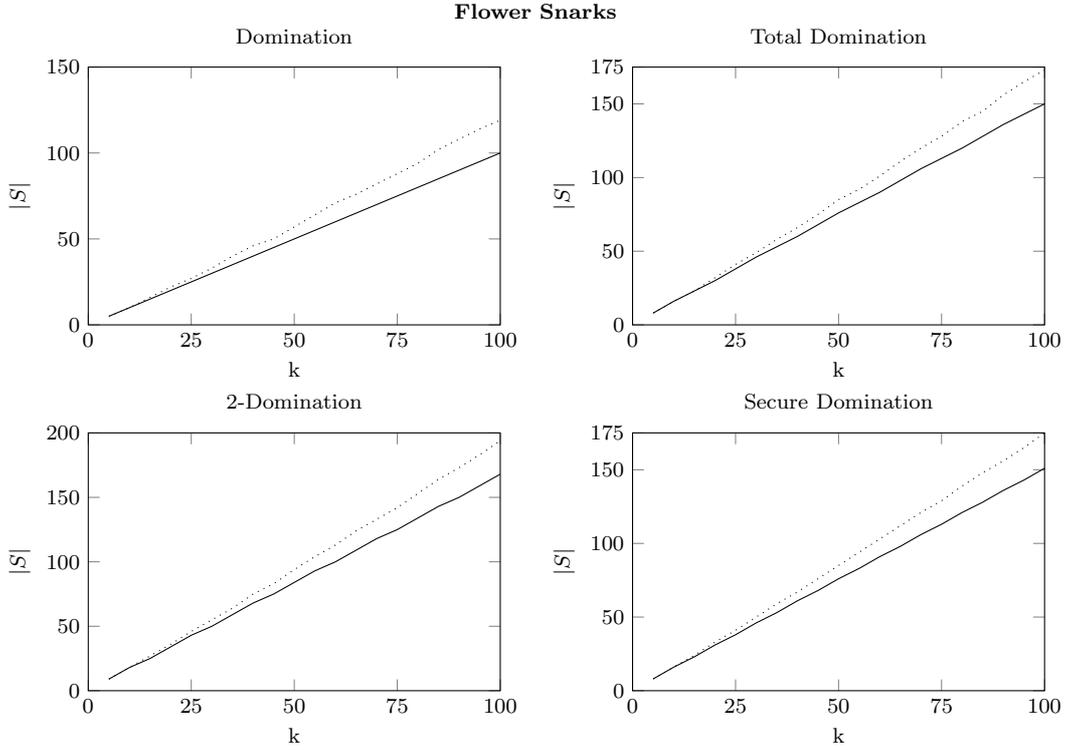

For square grid graphs and flower snarks, we present the results as plots since these graphs are not randomised, and follow a set structure. In Figures \ref{grid_plots}--\ref{flower_plots}, the dotted lines represent the best obtained values from the cross entropy method, and the solid lines represent the known optimal values. If the optimal values are not known, a dashed line is used to indicate the upper bound obtained by CPLEX after 10,000 seconds.

For the remaining graphs, results are displayed in Tables \ref{udg_table}--\ref{dimacs_table}. The ``CE" column lists the best obtained values from the cross entropy method, while the ``Sol" column lists either the optimal value (if it is known) or else the upper bound obtained by CPLEX after 10,000 seconds. In the latter case, the number is presented with an overline. Finally, the ``Gap" column gives the difference between these two numbers. If CPLEX has not even been able to produce an upper bound within 10,000 seconds, we indicate this with a dash (-).

%grids plots
\begin{figure}[h!]
    \centering
            {\bf \vspace*{2cm}Square Grid Graphs}\\
\begin{tikzpicture}
\begin{axis}[
title={Domination}, xlabel={n}, ylabel={$|S|$}, xmin=0, xmax=20,ymin=0, ymax=175, xtick={0, 5, 10, 15, 20}, ytick={0, 50, 100, 150}, ymode=linear, grid style=dashed,
]

\addplot[dotted]
coordinates {
(3,3)(4,4)(5,7)(6,10)(7,12)(8,16)(9,21)(10,26)(11,32)(12,38)(13,45)(14,52)(15,62)(16,69)(17,80)(18,89)(19,99)(20,110)
};

\addplot[solid]
coordinates {
(3,3)(4,4)(5,7)(6,10)(7,12)(8,16)(9,20)(10,24)(11,29)(12,35)(13,40)(14,47)(15,53)(16,60)(17,68)(18,76)(19,84)(20,92)
};

\end{axis}
\end{tikzpicture}
    \hskip 5pt
\begin{tikzpicture}
\begin{axis}[
title={Total Domination}, xlabel={n}, ylabel={$|S|$}, xmin=0, xmax=20,ymin=0, ymax=150, xtick={0, 5, 10, 15, 20}, ytick={0, 50, 100, 150}, ymode=linear, grid style=dashed,
]

\addplot[dotted]
coordinates {
(3,3)(4,6)(5,9)(6,12)(7,15)(8,20)(9,25)(10,31)(11,38)(12,45)(13,54)(14,62)(15,72)(16,82)(17,93)(18,104)(19,116)(20,129)
};

\addplot[solid]
coordinates {
(3,3)(4,6)(5,9)(6,12)(7,15)(8,20)(9,25)(10,30)(11,35)(12,42)(13,49)(14,56)(15,63)(16,72)(17,81)(18,90)(19,99)(20,110)
};

\end{axis}
\end{tikzpicture}
    \vspace*{0.5cm}
\begin{tikzpicture}
\begin{axis}[
title={2-Domination}, xlabel={n}, ylabel={$|S|$}, xmin=0, xmax=20,ymin=0, ymax=200, xtick={0, 5, 10, 15, 20}, ytick={0, 50, 100, 150, 200}, ymode=linear, grid style=dashed,
]

\addplot[dotted]
coordinates {
(3,4)(4,8)(5,11)(6,16)(7,21)(8,29)(9,36)(10,44)(11,54)(12,64)(13,76)(14,88)(15,101)(16,116)(17,131)(18,148)(19,164)(20,182)
};

\addplot[solid]
coordinates {
(3,4)(4,6)(5,11)(6,16)(7,21)(8,27)(9,34)(10,42)(11,50)(12,59)(13,69)(14,79)(15,90)(16,102)(17,114)(18,127)(19,141)(20,155)
};

\end{axis}
\end{tikzpicture}
    \hskip 5pt
\begin{tikzpicture}
\begin{axis}[
title={Secure Domination}, xlabel={n}, ylabel={$|S|$}, xmin=0, xmax=20,ymin=0, ymax=175, xtick={0, 5, 10, 15, 20}, ytick={0, 50, 100, 150}, ymode=linear, grid style=dashed,
]

\addplot[dotted]
coordinates {
(3,4)(4,7)(5,9)(6,13)(7,18)(8,24)(9,30)(10,38)(11,46)(12,55)(13,65)(14,76)(15,88)(16,100)(17,113)(18,126)(19,141)(20,157)
};

\addplot[solid]
coordinates {
(3,4.00)(4,7.00)(5,9.00)(6,13.00)(7,18.00)(8,23.00)(9,28.00)(10,35.00)(11,42.00)(12,49.00)
};

\addplot[dashed]
coordinates {
(12,49.00)(14,66.00)(15,77.00)(18,111.00)(19,123.00)(20,137.00)
};

\end{axis}
\end{tikzpicture}
    \caption{Results for square grid graphs $G(n,n)$. The size of the best solutions returned by the cross entropy method are displayed using a dotted line, while the known optimal values are displayed using a solid line. In the case of secure domination, a dashed line is used for the upper bounds obtained by CPLEX after 10,000 seconds.}
    %In the case of secure domination, the results are not known from literature, so a solid line is used for the values obtained by CPLEX, and a dashed line is used for the upper bounds obtained by CPLEX after 10,000 seconds.}
    \label{grid_plots}
\end{figure}
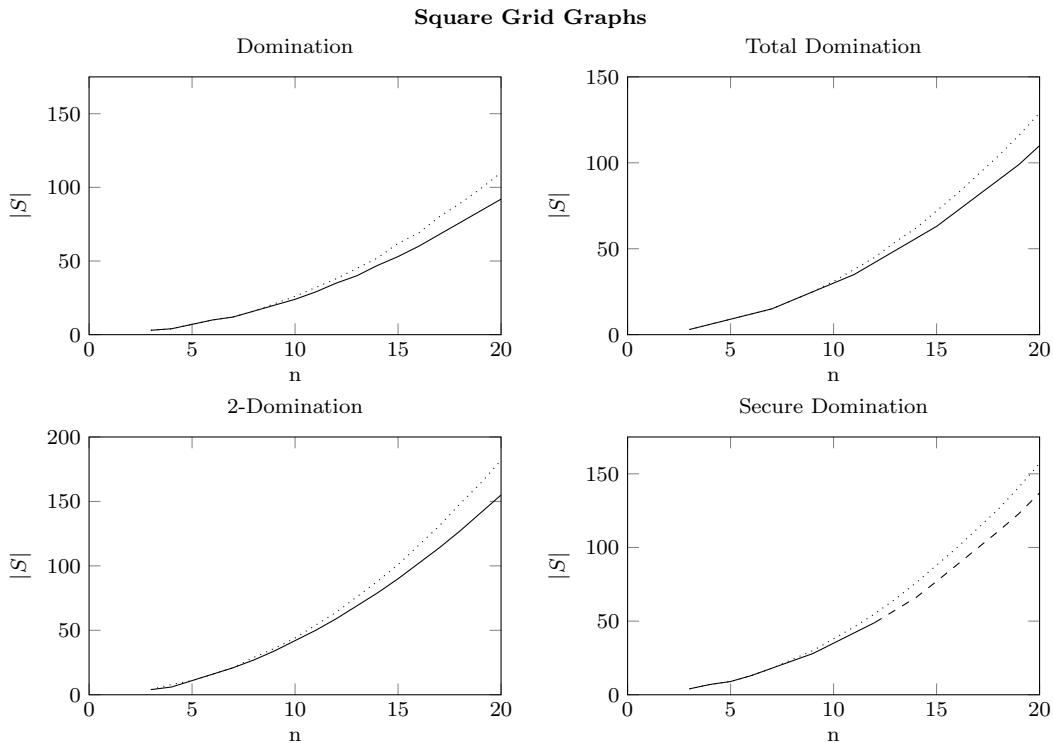

%flowers plots

%udg table
\begin{table}[h!]
\centering
\vspace*{1.5cm}\begin{tabular}{|l|ccc|ccc|ccc|ccc|}
\hline
\multirow{2}{*}{Instance} & \multicolumn{3}{c|}{Domination} & \multicolumn{3}{c|}{Total Domination} & \multicolumn{3}{c|}{2-Domination} & \multicolumn{3}{c|}{Secure Domination} \\
\cline{2-13}
& CE & Sol & Gap & CE & Sol & Gap & CE & Sol & Gap & CE & Sol & Gap \\
\hline
UDG\_100-0.7-10-10\_17 & 19 & 19 & 0 & 24 & 24 & 0 & 35 & 34 & 1 & 27 & 26 & 1 \\
UDG\_100-0.7-10-10\_32 & 18 & 18 & 0 & 25 & 25 & 0 & 35 & 34 & 1 & 27 & 27 & 0 \\
UDG\_100-0.7-10-10\_43 & 20 & 20 & 0 & 28 & 27 & 1 & 38 & 36 & 2 & 28 & 27 & 1 \\
UDG\_100-0.7-10-10\_55 & 20 & 20 & 0 & 24 & 24 & 0 & 36 & 36 & 0 & 27 & 27 & 0 \\
UDG\_100-0.7-10-10\_73 & 20 & 20 & 0 & 27 & 26 & 1 & 38 & 36 & 2 & 28 & 27 & 1 \\
UDG\_500-0.4-10-10\_0 & 65 & 55 & 10 & 85 & 71 & 14 & 125 & 106 & 19 & 94 & 81 & 13 \\
UDG\_500-0.4-10-10\_1 & 66 & 55 & 11 & 85 & 72 & 13 & 124 & 104 & 20 & 96 & 82 & 14 \\
UDG\_500-0.5-10-10\_0 & 44 & 37 & 7 & 58 & 46 & 12 & 85 & 71 & 14 & 67 & \textoverline{58} & 9 \\
UDG\_500-0.5-10-10\_1 & 46 & 37 & 9 & 60 & 47 & 13 & 87 & 73 & 14 & 68 & \textoverline{58} & 10 \\
UDG\_800-0.3-10-10\_0 & 118 & 99 & 19 & 154 & 124 & 30 & 217 & 183 & 34 & 166 & \textoverline{141} & 25 \\
UDG\_800-0.3-10-10\_1 & 118 & 99 & 19 & 150 & 122 & 28 & 217 & 184 & 33 & 169 & \textoverline{144} & 25 \\
UDG\_800-0.5-10-10\_0 & 48 & 39 & 9 & 63 & \textoverline{47} & 16 & 94 & \textoverline{74} & 20 & 74 & \textoverline{69} & 5 \\
UDG\_800-0.5-10-10\_1 & 49 & 39 & 10 & 63 & \textoverline{48} & 15 & 95 & \textoverline{76} & 19 & 75 & \textoverline{69} & 6 \\
UDG\_1000-0.3-10-10\_0 & 124 & 99 & 25 & 159 & \textoverline{123} & 36 & 227 & \textoverline{191} & 36 & 176 & \textoverline{149} & 27 \\
UDG\_1000-0.3-10-10\_1 & 121 & 99 & 22 & 155 & \textoverline{121} & 34 & 229 & \textoverline{200} & 29 & 177 & \textoverline{152} & 25 \\
UDG\_1000-0.5-10-10\_0 & 50 & 39 & 11 & 64 & \textoverline{51} & 13 & 97 & \textoverline{77} & 20 & 76 & \textoverline{75} & 1 \\
UDG\_1000-0.5-10-10\_1 & 52 & 39 & 13 & 64 & \textoverline{50} & 14 & 98 & \textoverline{76} & 22 & 77 & \textoverline{71} & 6 \\
\hline
\end{tabular}
\caption{Experimental results of the cross entropy method for Unit Disk graphs for each domination variant. The CE column shows the size of the best solutions returned by the cross entropy method. The Sol column shows the cardinality of an optimal set for the given instance, or an upper bound if the value is overlined. The Gap column shows the difference between the two values.}
\label{udg_table}
\end{table}

%random table
\begin{table}[h!]
\centering
\begin{tabular}{|l|ccc|ccc|ccc|}
\hline
\multirow{2}{*}{Instance} & \multicolumn{3}{c|}{Domination} & \multicolumn{3}{c|}{2-Domination} & \multicolumn{3}{c|}{Secure Domination} \\
\cline{2-10}
& CE & Sol & Gap & CE & Sol & Gap & CE & Sol & Gap \\
\hline
random100\_3 & 35 & 35 & 0 & 58 & 58 & 0 & 48 & 47 & 1 \\
random100\_4 & 29 & 28 & 1 & 48 & 48 & 0 & 41 & 39 & 2 \\
random100\_5 & 24 & 23 & 1 & 41 & 40 & 1 & 36 & 33 & 3 \\
random100\_6 & 21 & 19 & 2 & 35 & 33 & 2 & 31 & 28 & 3 \\
random250\_3 & 87 & 81 & 6 & 138 & 131 & 7 & 115 & 105 & 10 \\
random250\_4 & 73 & 65 & 8 & 117 & 108 & 9 & 102 & 90 & 12 \\
random250\_5 & 62 & 53 & 9 & 103 & 92 & 11 & 91 & - & - \\
random250\_6 & 56 & 46 & 10 & 92 & 81 & 11 & 82 & \textoverline{71} & 11 \\
random500\_3 & 168 & 148 & 20 & 274 & 255 & 19 & 234 & 209 & 25 \\
random500\_4 & 141 & 117 & 24 & 228 & 203 & 25 & 202 & \textoverline{170} & 32 \\
random500\_5 & 121 & 98 & 23 & 200 & \textoverline{171} & 29 & 179 & \textoverline{147} & 32 \\
random500\_6 & 108 & 84 & 24 & 180 & \textoverline{151} & 29 & 162 & \textoverline{134} & 28 \\
random800\_3 & 280 & 250 & 30 & 450 & 419 & 31 & 387 & 344 & 43 \\
random800\_4 & 234 & 195 & 39 & 383 & 339 & 44 & 336 & \textoverline{282} & 54 \\
random800\_5 & 207 & 163 & 44 & 339 & \textoverline{290} & 49 & 299 & \textoverline{249} & 50 \\
random800\_6 & 184 & \textoverline{141} & 43 & 301 & \textoverline{250} & 51 & 272 & \textoverline{222} & 50 \\
random1000\_3 & 374 & 329 & 45 & 578 & 536 & 42 & 504 & \textoverline{449} & 55 \\
random1000\_4 & 310 & 256 & 54 & 495 & 438 & 57 & 434 & \textoverline{367} & 67 \\
random1000\_5 & 269 & 212 & 57 & 430 & \textoverline{366} & 64 & 382 & \textoverline{321} & 61 \\
random1000\_6 & 242 & \textoverline{183} & 59 & 391 & \textoverline{323} & 68 & 351 & \textoverline{285} & 66 \\
\hline
\end{tabular}
\caption{Experimental results of the cross entropy method for Erd\H{o}s-R\'{e}nyi random graphs for each domination variant except for Total Domination (since these graphs usually contain isolated vertices). The CE column shows the size of the best solutions returned by the cross entropy method. The Sol column shows the cardinality of an optimal set for the given instance, or an upper bound if the value is overlined. The Gap column shows the difference between the two values. A dash is used when CPLEX was unable to even obtain an upper bound within 10,000 seconds.}
\label{random_table}
\end{table}

%literature instance table
\begin{table}[h!]
\centering
\vspace*{1.0cm}\begin{tabular}{|l|ccc|ccc|ccc|ccc|}
\hline
\multirow{2}{*}{Instance} & \multicolumn{3}{c|}{Domination} & \multicolumn{3}{c|}{Total Domination} & \multicolumn{3}{c|}{2-Domination} & \multicolumn{3}{c|}{Secure Domination} \\
\cline{2-13}
& CE & Sol & Gap & CE & Sol & Gap & CE & Sol & Gap & CE & Sol & Gap \\
\hline
adjnoun \cite{newman2006finding} & 18 & 18 & 0 & 19 & 19 & 0 & 39 & 38 & 1 & 32 & 31 & 1 \\
anna \cite{johnson1996cliques} & 12 & 12 & 0 & 12 & 12 & 0 & 47 & 47 & 0 & 42 & 42 & 0 \\
david \cite{johnson1996cliques} & 2 & 2 & 0 & 2 & 2 & 0 & 26 & 26 & 0 & 24 & 24 & 0 \\
dolphins \cite{lusseau2003bottlenose} & 14 & 14 & 0 & 17 & 17 & 0 & 27 & 27 & 0 & 22 & 22 & 0 \\
football \cite{girvan2002community} & 13 & 12 & 1 & 15 & 13 & 2 & 24 & 21 & 3 & 18 & 17 & 1 \\
gplus\_2000 \cite{chalupa2018order} & 236 & 170 & 66 & 188 & 181 & 7 & 1036 & 965 & 71 & 949 & - & - \\
gplus\_500 \cite{chalupa2018order} & 42 & 42 & 0 & 45 & 45 & 0 & 303 & 297 & 6 & 274 & 267 & 7 \\
homer \cite{johnson1996cliques} & 97 & 96 & 1 & $\infty$ & $\infty$ & & 323 & 317 & 6 & 294 & 282 & 12 \\
huck \cite{johnson1996cliques} & 9 & 9 & 0 & 11 & 11 & 0 & 21 & 21 & 0 & 15 & 15 & 0 \\
lesmis \cite{knuth1993stanford} & 10 & 10 & 0 & 10 & 10 & 0 & 33 & 33 & 0 & 28 & 28 & 0 \\
netscience \cite{newman2006finding} & 509 & 477 & 32 & $\infty$ & $\infty$ & & 933 & 915 & 18 & 643 & 623 & 20 \\
pokec\_2000 \cite{chalupa2018order} & 78 & 75 & 3 & 75 & 75 & 0 & 879 & 816 & 63 & 853 & - & - \\
pokec\_500 \cite{chalupa2018order} & 16 & 16 & 0 & 16 & 16 & 0 & 266 & 264 & 2 & 257 & 251 & 6 \\
polbooks \cite{chalupa2018order} & 14 & 13 & 1 & 15 & 15 & 0 & 24 & 22 & 2 & 21 & 19 & 2 \\
power \cite{watts1998collective} & 1747 & 1481 & 266 & 1932 & 1801 & 131 & 3002 & 2795 & 207 & 2575 & - & - \\
zachary \cite{zachary1977information} & 4 & 4 & 0 & 4 & 4 & 0 & 12 & 12 & 0 & 9 & 9 & 0 \\
\hline
\end{tabular}
\caption{Experimental results of the cross entropy method for the selected literature instances for each domination variant. The CE column shows the size of the best solutions returned by the cross entropy method. The Sol column shows the cardinality of an optimal set for the given instance, or an upper bound if the value is overlined. The Gap column shows the difference between the two values. A dash is used when CPLEX was unable to even obtain an upper bound within 10,000 seconds. An $\infty$ symbol is used whenever the graph contains no total dominating set.}
\label{dimacs_table}
\end{table}

%\cite{chalupa2018order,johnson1996cliques,newman2006finding,lusseau2003bottlenose,girvan2002community,knuth1993stanford,watts1998collective,zachary1977information}

%talk about results here i guess
\newpage

Across all of the experiments, it appears that the cross entropy algorithm described in Section \ref{sec-ce} produces solutions of similar quality; that is, for all kinds of graphs and variants of domination considered, it was generally able to find solutions within 10 to 25 percent of the known optimal solution or upper bound. However, certain variants of domination did occasionally perform slightly better for different graphs. For instance, secure domination performed relatively well for unit disk graphs, 2-domination performed relatively well for \ER random graphs, and total domination performed relatively well for the selected instances from the literature.

We highlight a few instances in particular. From Table \ref{dimacs_table}, the gplus\_2000 instance provided fascinating results. For (standard) domination, the best found solution had a gap of 66 from the optimal solution, and yet for total domination, the gap is only 7. Even for 2-domination, the gap is 71 but the size of the optimal 2-dominating set is more than five times as large as for the optimal dominating set. Similarly, for pokec\_2000, for (standard) domination the best found solution had a gap of 3, but for total domination an optimal solution was found which itself is also an optimal dominating set. It seems that in these more difficult instances, having a variant of domination can actually be beneficial; one possible explanation is that the more restrictive constraints of these variants force the cross entropy algorithm to consider only certain types of dominating sets which are perhaps closer to the optimal solution on average.

\begin{figure}[h!]
\centering
\begin{tikzpicture}
\begin{axis}[
title={Square Grid Times}, xlabel={n}, ylabel={time (s)}, xmin=3, xmax=20,ymax=10, xtick={5, 10, 15, 20}, ytick={0.0001, 0.001, 0.01, 0.1, 1, 10}, ymode=log, grid style=dashed, width=12cm, height=9cm
]

\addplot[solid]
coordinates {
(3,0.003)(4,0.009)(5,0.014)(6,0.02)(7,0.035)(8,0.053)(9,0.071)(10,0.098)(11,0.127)(12,0.209)(13,0.259)(14,0.293)(15,0.375)(16,0.413)(17,0.534)(18,0.825)(19,0.969)(20,1.055)
};

\addplot[dashed]
coordinates {
(3,0.006)(4,0.009)(5,0.013)(6,0.021)(7,0.041)(8,0.05)(9,0.076)(10,0.105)(11,0.151)(12,0.199)(13,0.24)(14,0.295)(15,0.412)(16,0.501)(17,0.523)(18,0.691)(19,0.864)(20,1.101)
};

\addplot[dotted]
coordinates {
(3,0.003)(4,0.009)(5,0.016)(6,0.028)(7,0.04)(8,0.052)(9,0.077)(10,0.123)(11,0.138)(12,0.161)(13,0.241)(14,0.309)(15,0.362)(16,0.474)(17,0.651)(18,0.818)(19,0.963)(20,1.1)
};

\addplot[dashdotted]
coordinates {
(3,0.023)(4,0.05)(5,0.109)(6,0.19)(7,0.279)(8,0.395)(9,0.566)(10,0.713)(11,1.047)(12,1.162)(13,1.44)(14,1.696)(15,1.874)(16,2.655)(17,3.181)(18,3.843)(19,5.147)(20,5.632)
};

\end{axis}

\end{tikzpicture}
\caption{Average run time in seconds for a single run of the cross entropy algorithm on square grid graphs $G(n,n)$, for the different variants of domination. The solid line corresponds to the domination problem; the dashed line to total domination; the dotted line to 2-domination; the dashdotted line to secure domination.}
\label{time_plot}
\end{figure}
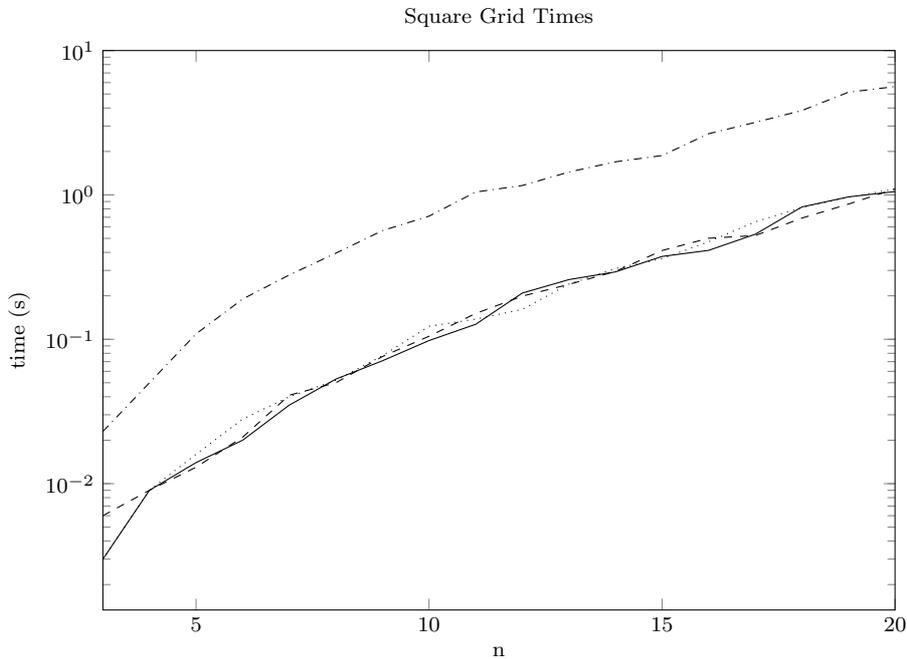

In both the gplus\_2000 and pokec\_2000 instances, CPLEX was unable to produce even an upper bound for secure domination within 10,000 seconds. This highlights the value of having a heuristic in hand which can be easily implemented for any sensible variant of domination and can generate solutions of reasonable quality in quick time. Across our experiments, we found that the size of the graph was the main factor in how long the cross entropy algorithm took to run, with the structure of the graph relatively unimportant. As such, in Figure \ref{time_plot} we provide average run times just for the square grid graphs, with very similar performance being exhibited for all other graphs.

%grid time plot

As can be seen, the time taken for the cross entropy algorithm to run is fairly consistent between domination, 2-domination and total domination, while it is roughly an order of magnitude slower for secure domination. This demonstrates that the main factor in the overall computation time for our implementation is just the time taken to actually generate the sets each iteration, with the remainder of the algorithm adding very little overhead per iteration. As such, it is imperative that the generation algorithms are optimised to ensure the overall algorithm can be run efficiently. We note here that even for graphs with 400 vertices, the cross entropy algorithm was able to produce dominating sets, 2-dominating sets and total dominating sets of reasonable quality in less than a second, and secure dominating sets in less then ten seconds, demonstrating the efficiency of this algorithm. Indeed, by looking at Figure \ref{grid_plots} we see that CPLEX was unable to find an optimal secure dominating set for the square grid graph $G(13,13)$, and after 10,000 seconds of computation it was only about to find an upper bound marginally better than that which was produced by the cross entropy implemented in less than 2 seconds.

\section{Conclusion and Future Work}\label{sec-conclusions}

Although there is a wealth of fast heuristics for the (standard) domination problem, its variants often have very few such heuristics available, and developing such heuristics typically requires one to first analyse the specific variant being considered. In contrast, the implementation of the cross entropy method described in this paper can be easily applied to most variants of domination without requiring any modification, and the experiments in Section \ref{sec-results} indicate that solution quality remains consistent for different variants. In particular, our method is suitable for any variant of domination in which simply adding vertices randomly is guaranteed to result in a set meeting the criteria of the variant; or, if there are no such sets in a given graph, we should be able efficiently identify this in advance. Then all that is required to use our method is to provide a {\em checking algorithm} which determines whether any given set meets the criteria of the desired variant.

The implementation itself is very lightweight outside of the checking algorithms, making it particularly suitable for large instances. Although it does not typically return an optimal solution for large instances, it is able to find good solutions very efficiently without the associated memory and performance issues that come with formulating the problem as a mixed-integer linear program, or requiring any other analysis of the underlying problem. The checking algorithm itself may be sped up in a number of ways; for example, the domination criteria do not need to be checked after each vertex is added and can instead be checked at regular intervals, particularly given that Phase 2 of Algorithm \ref{dom_alg} will trim unnecessary vertices out of the set. The checking algorithm may also be sped up using specialised hardware such as GPU to distribute the workload in a parallel sense. Since the checking algorithm is, by far, the slowest component of the algorithm, any improvements to it correspond to a direct improvement in the overall computation time.

We finish by noting that there are some variants of domination for which it is not the case that adding vertices randomly is guaranteed to result in a valid solution. For example, an {\em independent dominating set} is a dominating set which is also an independent set in the underlying graph, and adding vertices randomly could potentially violate the latter criterion. However, even in such cases, it may be possible to modify the cross entropy implementation accordingly. For example, in the case of independent domination, upon adding a vertex to $S$, in addition to setting the probability of that vertex to 0, one could also set the probability of its neighbours to 0, and this would ensure that an independent dominating set is generated.

\section*{Data Availability Statement}

The set of instances analysed in this study, and the results of that analysis, are available in the Github repository, \url{https://github.com/flinders-maths/test_instances}.

%recap the algorithms and its reason for existing
%go over its performance
%reccos for future work

%\begin{thebibliography}{9}
%\bibitem{adel} M. Guo, J. Li, A. Neumann, F. Neumann, H. Nguyen. Practical and Fixed-Parameter Algorithms for Defending Active Directory Style Attack Graphs, arXiv preprint arXiv:2112.13175 (2021).

% \end{thebibliography}

\bibliography{references.bib}
\bibliographystyle{acm}

\end{document}